\documentclass[12pt]{elsarticle}
\usepackage{amsfonts}
\usepackage{latexsym}
\usepackage{color}
\usepackage{graphicx}

\newtheorem{theorem}{Theorem}
\newtheorem{corollary}{Corollary}

\newtheorem{example}{Example}

\newtheorem{definition}{Definition}
\newtheorem{remark}{Remark}

\newcommand{\bfe}{\mbox{$\mbox{\boldmath $e$}$}} 
\newcommand{\bfg}{\mbox{$\mbox{\boldmath $g$}$}} 
\newcommand{\bfx}{\mbox{$\mbox{\boldmath $x$}$}} 
\newcommand{\bfy}{\mbox{$\mbox{\boldmath $y$}$}} 

\newcommand{\sbfx}{\mbox{\boldmath \footnotesize $x$}} 

\newcommand{\bnab}{\mbox{$\mbox{\boldmath $\nabla$}$}} 
\newcommand{\bxi}{\mbox{$\mbox{\boldmath $\xi$}$}} 

\newcommand{\bphi}{\mbox{$\mbox{\boldmath $\varphi$}$}} 
\newcommand{\bpi}{\mbox{$\mbox{\boldmath $\pi$}$}} 

\newcommand{\bfzero}{\mbox{$\mbox{\boldmath $0$}$}} 
\newcommand{\bfone}{\mbox{$\mbox{\boldmath $1$}$}} 

\date{}

\journal{Journal of Mathematical Analysis and Applications}

\begin{document}

\begin{frontmatter}

\title{The fundamental equations for inversion of operator pencils on Banach space}


\author[scg]{Amie Albrecht}
\author[scg]{Phil Howlett\corref{cor1}}
\author[sms]{Charles Pearce\corref{cor2}}

\address[scg]{Scheduling and Control Group, Centre for Industrial and Applied Mathematics, University of South Australia, Mawson Lakes 5095, Australia}
\address[sms]{School of Mathematical Sciences, University of Adelaide, Adelaide 5000, Australia}
\cortext[cor1] {Corresponding author. Email: phil.howlett@unisa.edu.au}
\cortext[cor2]{This paper is dedicated to our co--author, Charles E. M. Pearce, who died on June 8th 2012, in a tragic car accident on the South Island of New Zealand.   Charles was a most generous and agreeable person.  He is sadly missed by family, friends and colleagues.}

\begin{abstract}
\label{abs}

\noindent We prove that the resolvent of a linear operator pencil is analytic on an open annulus if and only if the coefficients of the Laurent series satisfy a system of fundamental equations and are geometrically bounded.  Our analysis extends earlier work on the fundamental equations to include the case where the resolvent has an isolated essential singularity.  We find a closed form for the resolvent and use the fundamental equations to establish key spectral separation properties when the resolvent has only a finite number of isolated singularities.  Finally we show that our results can also be applied to polynomial pencils.

\end{abstract}

\begin{keyword}
operator pencil \sep resolvent \sep fundamental equations \sep singular perturbation
\MSC[2010] 47A10 \sep 47A55 \sep 47A56
\end{keyword}

\end{frontmatter}

\section{Introduction}
\label{intro}

Let $H,K$ be Banach spaces and let $A_0, A_1 \in {\mathcal L}(H,K)$, where $A_0$ is singular.  Define a linear operator pencil $A: {\mathbb C} \rightarrow {\mathcal L}(H,K)$ by the formula $A(z) = A_0 + A_1z$.  We show that the resolvent $R:\mathcal U_{s,r} \rightarrow {\mathcal L}(K,H)$ for $A$ on the annulus ${\mathcal U}_{s,r} = \{z \in {\mathbb C} \mid s < |z| < r\}$ is defined by a Laurent series $R(z) = \sum_{j \in {\mathbb Z}} R_jz^j$  if and only if the coefficients $R_j \in {\mathcal L}(K,H)$ for each $j \in {\mathbb Z}$ satisfy a system of left and right fundamental equations.  If a solution exists we show that the solution is uniquely defined by the coefficients $R_{-1}$ and $R_0$ and we find a closed form $R(z) = (I z + R_{-1}A_0)^{-1}R_{-1} + (I + R_0A_1z)^{-1}R_0$ for the resolvent.  It follows too that the functions ${\mathcal R}_{\lambda} = \lambda^{-1}R(-\lambda^{-1})A_0 \in {\mathcal L}(H)$ and ${\mathcal S}_{\lambda} = \lambda^{-1}A_0R(-\lambda^{-1}) \in {\mathcal L}(K)$ each satisfy a classical resolvent equation.  We use the fundamental equations to show that the operators $P = R_{-1}A_1 \in {\mathcal L}(H)$ and $Q = A_1R_{-1}\in {\mathcal L}(K)$ define corresponding projections that separate the bounded and unbounded components $\sigma_1 \subset \{ z \in {\mathbb C} \mid |z| \leq s \}$ and $\sigma_2 \subset \{ z \in {\mathbb C} \mid |z| \geq r \}$ of the spectral set $\sigma$ and we use these projections to establish useful properties of the coefficients $\{R_j\}_{j \in {\mathbb Z}}$.  We extend our main results to polynomial pencils and investigate the global structure of $R(z)$ in the case where $R(z)$ has only a finite number of isolated singularities.

\subsection{Motivation: Mean first-passage times for perturbed Markov processes}
\label{mpt}

Determination of mean first-passage times in perturbed Markov processes can be solved by a linear pencil inversion.  The intrinsic structure of a Markov process can be substantially changed by a small perturbation.  For instance the perturbation may introduce state transitions that are not possible in the unperturbed process.  For a Markov process $T$ it is known that the mean first-passage times between states can be calculated by finding the linear operator
\begin{equation}
\label{mpt1}
[I - T + T^{\infty}]^{-1}
\end{equation}
where $T^{\infty} = \lim_{n \rightarrow \infty} T^n$ is suitably defined.  We consider a finite-dimensional problem on a discrete state space and refer the reader to a previous paper \cite{how4} for an analogous infinite-dimensional problem on a continuous state space.  Let $T:{\mathbb R}^{1 \times (r+1)} \rightarrow {\mathbb R}^{1 \times (r+1)}$ be a Markov process defined by $T(\bpi) = \bpi P$ for each probability vector $\bpi \in {\mathbb R}^{1 \times (r+1)}$ where
\begin{equation}
\label{mpt2}
P = \left[\begin{array}{ccccc}
\rule{0cm}{0.5 cm} 1 & 0 &  \cdots & 0 & 0 \\
\rule{0cm}{0.5 cm} \frac{1}{2} & \frac{1}{2} & \cdots & 0 & 0 \\
\rule{0cm}{0.5 cm} \vdots & \vdots & \ddots & \vdots & \vdots \\
\rule{0cm}{0.5 cm} \frac{1}{r} & \frac{1}{r} & \cdots & \frac{1}{r} & 0 \\
\rule{0cm}{0.5 cm} \frac{1}{r+1} & \frac{1}{r+1} & \cdots & \frac{1}{r+1} & \frac{1}{r+1} \end{array} \right] \in {\mathbb R}^{(r+1) \times (r+1)}
\end{equation}
Now let $T_{\epsilon}:{\mathbb R}^{1 \times (r+1)} \rightarrow {\mathbb R}^{1 \times (r+1)}$ be a perturbed process defined by
\begin{equation}
\label{mpt3}
T_{\epsilon}(\bpi) = [(1 - \epsilon)I + \epsilon T ](\bpi)
\end{equation}
where $I \in {\mathbb R}^{(r+1) \times (r+1)}$ is the identity matrix and $0 \leq \epsilon \leq 1$.  The process $T_{\epsilon}$ is a perturbation of the identity.  It is a singular perturbation because the process changes radically for $\epsilon \neq 0$.  When $\epsilon = 0$ the Markov process is simply an identity transformation and the initial state does not change.  If we write $\bfe_1 = [1,0,\ldots,0] \in {\mathbb R}^{1 \times n}$ for each $n \in {\mathbb N}$ then it can be seen that $T_{\epsilon} (\bfe_1) = \bfe_1$ and hence, when $\epsilon > 0$, it follows that $T_{\epsilon}^n (\bpi) \rightarrow \bfe_1$ as $n \rightarrow \infty$.  If we regard the state space as the set of numbers
\begin{equation}
\label{mpt4}
S = \{ 0, \frac{1}{r}, \frac{2}{r}, \ldots, \frac{r-1}{r}, 1 \}
\end{equation}
then the invariant measure for the perturbed process $T_{\epsilon}$ lies entirely at zero.  The perturbed transformation $T_{\epsilon}$ allows leakage back to the zero state.  To find the mean first-passage times we must calculate
\begin{equation}
\label{mpt5}
[I - T_{\epsilon} + T_{\epsilon}^{\infty}]^{-1}.
\end{equation}
Define transformations $U_0, U_1:{\mathbb R}^{1 \times (r+1)} \rightarrow {\mathbb R}^{1 \times (r+1)}$ by setting $U_0(\bxi) = \bxi A_0^T$ and $U_1(\bxi) = \bxi A_1^T$ where $A_0^T = [\bfone^T,\bfzero^T,\ldots,\bfzero^T]  \in {\mathbb R}^{(r+1) \times (r+1)}$ and $A_1^T = I - P \in {\mathbb R}^{(r+1) \times (r+1)}$ and where we have used the notation $\bfzero = [0, \ldots, 0] \in {\mathbb R}^{1 \times (r+1)}$ and $\bfone = [1,\ldots,1] \in {\mathbb R}^{1 \times (r+1)}$.   Now $[I - T_{\epsilon} + T_{\epsilon}^{\infty}]^{-1} = [U_0 + \epsilon U_1]^{-1}$ and the required inversion\footnote{The operator is transposed at this stage in order to use the more conventional column vector notation in our proposed Gaussian elimination.} is equivalent to finding the resolvent operator for the matrix pencil $A(\epsilon) = A_0 + A_1 \epsilon$.  If we assume that the singularity at $\epsilon = 0$ is a first-order pole then it has been shown in \cite{how5} that for matrix problems the solution $\{R_j\}$ can be found by the Gaussian reduction
\begin{equation}
\label{mpt6}
\left[ \begin{array}{cccc|c}
A_0 & 0 & 0 & \cdots & 0 \\
A_1 & A_0 & 0 & \cdots & I \\
0 & A_1 & A_0 & \cdots & 0 \\
\vdots & \vdots & \vdots &  & \vdots \end{array} \right] \quad \longrightarrow \left[ \begin{array}{cccc|c}
I & 0 & 0 & \cdots & R_{-1} \\
0 & I & 0 & \cdots & R_0 \\
\vdots & \vdots & \vdots & & \vdots \end{array} \right].
\end{equation}
This yields the solution
\begin{equation}
\label{mpt7}
R_{-1} = \left[ \begin{array}{ccccc}
\rule{0cm}{0.5cm} 0 & -(s_1+1) & -(s_2+1)  & \cdots & -(s_r +1) \\
\rule{0cm}{0.5cm} 0 & 1+1 & 1 & \cdots & 1 \\
\rule{0cm}{0.5cm} 0 & 0 & 1+\frac{1}{2} & \cdots & \frac{1}{2} \\
\rule{0cm}{0.5cm} \vdots & \vdots & \vdots & \ddots & \vdots \\
\rule{0cm}{0.5cm} 0 & 0 & 0 & \cdots & 1+\frac{1}{r} \end{array} \right], \hspace{0.2cm} R_0 = \left[ \begin{array}{cccc}
1 & 1 & \cdots & 1 \\
0 & 0 & \cdots & 0 \\
\vdots & \vdots & \ddots & \vdots \\
0 & 0 & \cdots & 0 \end{array} \right],
\end{equation}
with $R_j = 0$ for $j \in {\mathbb N}$ where $s_n = \sum_{j=1}^n 1/j$ for each $n \in {\mathbb N}$.  The resolvent is given by
\begin{equation}
\label{mpt8}
(A_0 + A_1 \epsilon)^{-1} = \frac{R_{-1}}{\epsilon} + R_0
\end{equation}
for $0 < \epsilon \leq 1$.  Clearly $R(\epsilon) = (A_0 + A_1 \epsilon)^{-1}$ has a first-order pole at $\epsilon = 0$. 

\subsection{Motivation: Input retrieval in linear control systems}
\label{ir}

The problem of input retrieval in linear control systems can be solved by inversion of a linear pencil.  See \cite{how1,sai1} for an extensive discussion of the finite-dimensional problem and \cite{how5} for an outline of the infinite-dimensional problem.

\subsection{Previous work on spectral theory for linear operator pencils}
\label{pwst}

We outline relevant results from the paper by F. Stummel \cite{stu2}.   Stummel considers the spectral set $\sigma = \{ z \in {\mathbb C} \mid A(z)^{-1} \notin {\mathcal L}(H,K) \}$ for the linear pencil $A(z) = A_0 + A_1z \in {\mathcal L}(H,K)$ and uses line integrals to define complementary projections that isolate separate components of the spectrum.  In the case where $R(z) = A(z)^{-1} \in {\mathcal L}(K,H)$ is analytic on an annular region ${\mathcal U}_{s,r} = \{ z \in {\mathbb C} \mid s < |z| < r \}$ the spectral set has two disjoint components---an {\em interior} compact subset $\sigma_1 \subset \{ z \in {\mathbb C} \mid |z| \leq s \}$ and an {\em exterior} closed subset $\sigma_2 \subset \{ z \in {\mathbb C} \mid |z| \geq r \}$.  Following Stummel we define an operator $R_{-1} \in {\mathcal L}(K,H)$ by the contour integral formula
\begin{equation}
\label{pwst1}
R_{-1} = \frac{1}{2 \pi i} \int_{\Gamma_{\rho}} R(\zeta) d\zeta
\end{equation}
where $\Gamma_{\rho} = \{ \zeta = \rho e^{i \theta}\mid \theta \in [0, 2 \pi) \}$ and  $s < \rho < r$.  Stummel defines corresponding projections $P = R_{-1}A_1 \in {\mathcal L}(H)$ and $Q = A_1R_{-1} \in {\mathcal L}(K)$ to establish the formal separation of $\sigma_1$ and $\sigma_2$.  The method proposed by Stummel is described clearly by Gohberg {\em et al.} \cite[pp. 49--54]{goh1}.  The earlier underlying concept of spectral separation for a bounded linear operator is presented elegantly in the classic book by Kato \cite[pp. 178--179]{kat1}.  Although the integral formula (\ref{pwst1}) can be used to establish the existence of the key projection operators $P$ and $Q$ it cannot be used {\em per se} to calculate $R_{-1}$ because the integral depends on $R(\zeta) = \sum_{j \in {\mathbb Z}} R_j \zeta^j$ which in turn depends on $R_{-1}$.  

\subsection{Previous work on the fundamental equations}
\label{pwfe}

The use of a semi-infinite system of fundamental matrix equations to define a matrix power series inversion was suggested by Sain and Massey \cite{sai1} and later developed by Howlett \cite{how1} to solve the problem of input retrieval in finite-dimensional linear control systems.  Solution methods for the semi-infinite system of fundamental matrix equations were further developed by Avrachenkov and Lasserre \cite{avr1} to analyse singularly perturbed Markov chains and by Avrachenkov et al. \cite{avr2} to develop efficient numerical algorithms for the inversion of analytic matrix pencils.  For an analytic matrix pencil in the form
\begin{equation}
\label{pwfe1}
A(z) = A_0 + A_1z + A_2z^2 + \cdots
\end{equation}
on some open disk $|z| < r$, where $A_j \in {\mathbb C}^{n \times n}$, it is well known that if $X(z) = A(z)^{-1}$ has an isolated singularity at $z=0$ then the singularity must be a pole of finite order $s \leq n$.  Thus, in the notation of \cite{alb1,how5}, the inverse takes the form
\begin{equation}
\label{pwfe2}
X(z) = \frac{1}{z^s} \left[ X_0 + X_1z + X_2z^2 + \cdots \right]
\end{equation}
for some $s \in \{1,\ldots,n\}$ and hence, by equating coefficients in the equations $X(z)A(z) = A(z)X(z) = I$, one obtains a semi-infinite system of left and right fundamental matrix equations.  Thus the inverse $X(z)$ can conceivably be obtained by solving a system of linear matrix equations.  A key aspect of the solution, used by all authors \cite{avr1, avr2, how1, sai1, sch1}, is that rank considerations allow the semi-infinite system to be reduced to an equivalent finite system.  Although these considerations no longer apply in an infinite-dimensional setting the fundamental equations have nevertheless been used to define the inversion of operator pencils when the resolvent $R(z) = A(z)^{-1}$ has a finite-order pole at $z=0$.  For a first-order pole Howlett et al. \cite{how5} showed that a solution to the fundamental equations could be used to define the resolvent operator in terms of corresponding complementary projections in the domain and range spaces.  The same authors subsequently established analogous results for higher-order poles \cite{alb1}.  The primary purpose of this paper is to show that a solution to the fundamental equations is necessary and sufficient to define the resolvent---even when there is an essential singularity at the origin.  The previous papers \cite{alb1,how5} relied on the finite termination of the series of negative powers in the Laurent expansion to establish the existence of the key projection operators.  A key factor in the extension to an essential singularity where the series of negative powers does not terminate is the realisation that the bounds imposed by the inner and outer radii of convergence ensure that certain basic operator products must be zero.

\subsection{Other relevant work}

Recent work on the inversion of linearly-perturbed operators on Hilbert space by Howlett {\em et al.} \cite{how4} was inspired by the work of Schweitzer and Stewart \cite{sch1} on a corresponding matrix inversion problem.  Whereas Schweitzer and Stewart defined the resolvent using an algebraic separation of the underlying spaces into complementary subspaces the extension to operators on Hilbert space required a geometric separation using orthogonal projections.  Incidentally it is noted in passing that Howlett et al. \cite{how4} showed how spaces of Sobolev type could be used to extend the construction of the resolvent to unbounded operators.  Although separation of the underlying spaces as a direct sum is guaranteed in Hilbert space by the use of orthogonal projections this is no longer true in Banach space.  Thus, in Banach space, some other rationale is needed to define the required linear projections and establish the formal spectral separation.  The fundamental equations provide that rationale.  Finally we note that there is a substantial literature on the inversion of matrix pencils.  See \cite{avr1,wil1} for more information and additional references.

\section{The main results}
\label{mr}

If the series $R(z) = \sum_{j \in {\mathbb Z}} R_jz^j$ converges on the region ${\mathcal U}_{s,r}$ then for each $\delta, \epsilon$ with $s < s+ \epsilon < r-\delta < r$ there are constants $c_{\delta}, d_{\epsilon}$ such that
\begin{equation}
\label{gb}
\|R_j\| \leq c_{\delta} / (r - \delta)^j \quad \mbox{and} \quad \|R_{-j}\| \leq d_{\epsilon} (s+\epsilon)^j \quad
\mbox{for} \quad j \in {\mathbb N}.
\end{equation}
If we equate coefficients for each power of $z$ in the identities $R(z)A(z) = I \in {\mathcal L}(H)$ and $A(z)R(z) = I \in {\mathcal L}(K)$ then it follows that the sequence of coefficients $\{R_j\}_{j \in {\mathbb Z}} \in {\mathcal L}(K,H)$ is a solution to the left fundamental equations
\begin{equation}
\label{lfe}
R_{j-1}A_1 + R_jA_0 = \left \{ \begin{array}{ll}
I & \mbox{if}\ j=0 \\
0 & \mbox{if}\ j \neq 0 \end{array} \right.
\end{equation}
and the right fundamental equations
\begin{equation}
\label{rfe}
 A_1R_{j-1} + A_0R_j = \left \{ \begin{array}{ll}
I & \mbox{if}\ j=0 \\
0 & \mbox{if}\ j \neq 0. \end{array} \right.
\end{equation}
Conversely if the sequence of coefficients $\{R_j\}_{j \in {\mathbb Z}} \in {\mathcal L}(K,H)$ is a solution to the systems (\ref{lfe}) and (\ref{rfe}) and satisfies the bounds given in (\ref{gb}) then $R(z) = \sum_{j \in {\mathbb Z}} R_jz^j$ is well defined on ${\mathcal U}_{s,r}$ and satisfies the identities $R(z)A(z) = I$ and $A(z)R(z) = I$.  Hence we have established our first key result.

\begin{theorem}
\label{fundeqn}
Let $s,r \in {\mathbb R}$ with $0 \leq s < r$.  There exists an analytic resolvent $R:{\mathcal U}_{s,r} \rightarrow {\mathcal L}(K,H)$ for $A$ defined by $R(z) = \sum_{j \in {\mathbb Z}} R_j z^j$ if and only if $\{R_j\}_{j \in {\mathbb Z}} \in {\mathcal L}(K,H)$ are geometrically bounded by $(\ref{gb})$ and satisfy the fundamental equations $(\ref{lfe})$ and $(\ref{rfe})$.
\end{theorem}

\begin{corollary}
\label{fundeqnunique}
If $\{S_j\}_{j \in {\mathbb Z}}$ and $\{R_j\}_{j \in {\mathbb Z}}$ are solutions to the fundamental equations satisfying  $(\ref{gb})$, $(\ref{lfe})$ and $(\ref{rfe})$ then $S_j = R_j$ for all $j \in {\mathbb Z}$.
\end{corollary}

\noindent{\bf Proof.} \quad Let $R(z) = \sum_{j \in {\mathbb Z}} R_jz^j$ and $S(z) = \sum_{j \in {\mathbb Z}} S_jz^j$ for $z \in {\mathcal U}_{s,r}$.  By Theorem \ref{fundeqn} we have $S(z) = S(z)[A(z)R(z)] = [S(z)A(z)]R(z) = R(z)$.  $\hfill \Box$\medskip 

By considering the fundamental equations more closely we can establish the important properties of the coefficients $\{R_j\}_{j \in {\mathbb Z}}$.  We have our second key result\footnote{Note that condition $(i)$ of Theorem \ref{fundsol} is incorrectly stated in the published paper\textemdash Amie Albrecht, Phil Howlett and Charles Pearce (2014), The fundamental equations for inversion of operator pencils on Banach space, {\em Journal of Mathematical Analysis and Applications}, {\bf 413}, 411-421.}.

\begin{theorem}
\label{fundsol}
The coefficients $\{R_j\}_{j \in {\mathbb Z}} \in {\mathcal L}(K,H)$ satisfy $(\ref{gb})$, $(\ref{lfe})$ and $(\ref{rfe})$ if and only if the following are all satisfied: $(i)$ $P = R_{-1}A_1 \in {\mathcal L}(H)$ and $I - P = R_0A_0 \in {\mathcal L}(H)$ are complementary projections on $H$; and $Q = A_1R_{-1} \in {\mathcal L}(K)$ and $I - Q = A_0R_0 \in {\mathcal L}(K)$ are corresponding complementary projections on $K$;  $(ii)$ $A_i = QA_iP + (I-Q)A_i(I-P)$ for $i=0,1$;  $(iii)$ $R_j = PR_jQ$ for $j \leq -1$ and $R_j = (I-P)R_j(I-Q)$ for $j \geq 0$; $(iv)$ $R_{-j} = (-1)^{j-1}(R_{-1}A_0)^{j-1}R_{-1}$ and $R_j = (-1)^j (R_0A_1)^jR_0$ for $j \in {\mathbb N}$; and $(v)$ $\lim_{j \rightarrow \infty} \|(R_{-1}A_0)^j\|^{1/j} \leq s$  and $\lim_{j \rightarrow \infty} \|(R_0A_1)^j \|^{1/j} \leq 1/r$.
\end{theorem}

\noindent {\bf Proof.} \quad  Suppose $\{R_j\}_{j \in {\mathbb Z}} \in {\mathcal L}(K,H)$ satisfies $(\ref{gb})$, $(\ref{lfe})$ and $(\ref{rfe})$.  Then for $k>0$ and $\ell \geq 0$ the fundamental equations show that $R_{-k}A_0R_{\ell}=-R_{-(k+1)}A_1R_{\ell}=R_{-(k+1)}A_0R_{\ell +1}$ and so, by repeated application of this process, it follows that $R_{-k}A_0R_{\ell}=R_{-(k+j)}A_0R_{\ell +j}$ for all $j \in {\mathbb N}$.  Since
\begin{equation}
\label{fs1}
\|R_{-(k+j)}A_0R_{\ell+j}\| \leq \|A_0\| c_{\delta}d_{\epsilon} (s+\epsilon)^k(r-\delta)^{-{\ell}} \left( (s+\epsilon)/(r-\delta) \right)^j \rightarrow 0
\end{equation}
as $j \rightarrow \infty$ it follows that $R_{-k}A_0R_{\ell}=0$ for $k>0$ and $\ell \geq 0$.  Consequently $(R_0A_0-I)R_{\ell}=-R_{-1}A_1R_{\ell}=R_{-1}A_0R_{\ell +1}=0$ and hence $R_0A_0R_{\ell} = R_{\ell}$ for $\ell \geq 0$.  Similar arguments show that $R_{\ell}A_0R_{-k} = 0$ for $k>0$ and $\ell \geq 0$ and that $R_{\ell}A_0R_0 = R_{\ell}$ for $\ell \geq 0$.  These results can now be used directly to show that $R_{-k}A_1R_{\ell} = 0$ and $R_{\ell}A_1R_{-k} = 0$ for $k > 0$ and $\ell \geq 0$.   

To prove $(i)$, define $P = R_{-1}A_1 \in {\mathcal L}(H)$ and $Q = A_1R_{-1} \in {\mathcal L}(K)$.  Now $P - P^2 = P(I - P)  =  R_{-1}A_1R_0A_0 = 0$, so $P$ is a projection.  A similar argument establishes that $Q$ is a projection.  Hence $(i)$ is true. 

To prove $(ii)$, observe that $QA_i(I-P) = A_1R_{-1}A_iR_0A_0 = 0$ and $(I-Q)A_iP = A_0R_0A_iR_{-1}A_1 = 0$ for $i=0,1$.  Therefore $A_i = [Q + (I - Q)]A_i[P + (I - P)] = QA_iP + (I-Q)A_i(I-P)$. 

For $(iii)$ choose $k>0$ and $\ell \geq 0$.  We have $PR_{-k}Q =  (I - R_0A_0)R_{-k}(I-A_0R_0) = R_{-k}(I-A_0R_0) = R_{-k}$ and $(I-P)R_{\ell}(I-Q) =  R_0A_0R_{\ell}A_0R_0 = R_0A_0R_{\ell} =R_{\ell}$.  

To prove $(iv)$, observe that $R_j = R_0A_0R_j = (-1)R_0A_1R_{j-1}$ for each $j \in {\mathbb N}$ and that repeated application of this formula gives $R_j = (-1)^j(R_0A_1)^j R_0$ for $j \geq 0$.  A similar argument shows that $R_{-j} = (-1)^{j-1}(R_{-1}A_0)^{j-1}R_{-1}$ for each $j \in {\mathbb N}$.  

For $(v)$ with sufficiently small $\epsilon > 0$, we have
$$
\|(R_{-1}A_0)^n\|^{1/n} = \|R_{-n}A_0\|^{1/n} \leq (\|A_0\| d_{\epsilon})^{1/n} (s + \epsilon)
$$
for all $n \in {\mathbb N}$ and so
$$
\lim_{n \rightarrow \infty} \|(R_{-1}A_0)^n\|^{1/n} \leq s+\epsilon.
$$
Since $\epsilon>0$ can be chosen arbitrarily small the desired result follows.  The other inequality is established similarly.

Conversely, suppose $(i)$--$(v)$ are all true.  From $(iv)$ and $(v)$ it follows that $R(z) = \sum_{j \in {\mathbb Z}} R_jz^j$
is well defined and analytic on ${\mathcal U}_{s,r}$.  From $(i)$ it can be seen that $R_{-1}A_1+ R_0A_0 = I$ and $A_1R_{-1} + A_0R_0 = I$.  When $k > 0$ and $\ell \geq 0$ the conditions $(ii)$--$(iv)$ imply
$$
R_{-k}A_iR_{\ell} = PR_{-k}Q \left[QA_iP + (I-Q)A_i(I-P) \right](I-P)R_{\ell}(I-Q) = 0
$$
and
$$
R_{\ell}A_iR_{-k} = (I-P)R_{\ell}(I-Q)\left[QA_iP + (I-Q)A_i(I-P) \right]PR_{-k}Q = 0
$$
for each $i=0,1$.  In particular we have $R_{-1}A_0R_0=0$ and $R_0A_1R_{-1}=0$.  Therefore for $k,\ell \in {\mathbb N}$ we have
$$
R_{-k-1}A_1 + R_{-k}A_0 = (-1)^{k-1}(R_{-1}A_0)^k[I - R_{-1}A_1] = (-1)^{k-1}(R_{-1}A_0)^kR_0A_0 = 0
$$
and
$$
R_{\ell-1}A_1 + R_{\ell}A_0 =  (-1)^{\ell-1}(R_0A_1)^{\ell}[I - R_0A_0]  = (-1)^{\ell-1}(R_0A_1)^{\ell-1}R_{-1}A_1 = 0.
$$
A similar argument provides $A_1R_{\ell-1} + A_0R_{\ell} = 0$ and $A_1R_{-k-1} + A_0R_{-k} = 0$.  Thus the fundamental equations (\ref{lfe}) and (\ref{rfe}) are satisfied.  Since the bounds (\ref{gb}) are also satisfied, it follows that $R(z) = A(z)^{-1}$ for $z \in {\mathcal U}_{s,r}$.  $\hfill \Box$

\begin{corollary}
\label{basicsol}
Let $s,r \in {\mathbb R}$ with $0 \leq s < r$.  The resolvent $R:{\mathcal U}_{s,r} \rightarrow {\mathcal L}(K,H)$
is analytic if and only if there exist $R_{-1}, R_0 \in {\mathcal L}(K,H)$ such that $(i)$ $R_{-1}A_1 + R_0A_0 = I$ and $A_1R_{-1} + A_0R_0 = I$;  $(ii)$ $R_{-1}A_iR_0 = 0$ and $R_0A_iR_{-1} = 0$ for each $i = 0,1$; and $(iii)$ $\lim_{j \rightarrow \infty}  \|(R_{-1}A_0)^j\|^{1/j} \leq s$ and $\lim_{j \rightarrow \infty} \|(R_0A_1)^j \|^{1/j} \leq 1/r$.  If these conditions are satisfied then $R_{-1}, R_0$ are uniquely defined and
\begin{eqnarray}
\label{cf}
R(z) & = & R_{[0]}(z) + R_{[0]}^c(z) = PR(z)Q + (I-P)R(z)(I-Q) \hspace {2.5cm} \nonumber \\
& & \hspace {2.5cm} = (I z + R_{-1}A_0)^{-1}R_{-1} + (I + R_0A_1z)^{-1}R_0
\end{eqnarray}
for all $z \in {\mathcal U}_{s,r}$.  If we define ${\mathcal R}_{\lambda} \in {\mathcal L}(H)$, ${\mathcal S}_{\lambda} \in {\mathcal L}(K)$ by ${\mathcal R}_{\lambda} = \lambda^{-1}R(-\lambda^{-1})A_0$ and ${\mathcal S}_{\lambda} = \lambda^{-1}A_0R(-\lambda^{-1})$, then ${\mathcal R}_{\lambda}, {\mathcal S}_{\lambda}$ satisfy the respective resolvent equations ${\mathcal R}_{\lambda} - {\mathcal R}_{\mu} = (\mu - \lambda) {\mathcal R}_{\lambda} {\mathcal R}_{\mu}$ and ${\mathcal S}_{\lambda} - {\mathcal S}_{\mu} = (\mu - \lambda) {\mathcal S}_{\lambda} {\mathcal S}_{\mu}$ for $\lambda, \mu \in {\mathcal U}_{r^{-1},s^{-1}}$.
\end{corollary}

\noindent {\bf Proof.}\quad If $R(z) = A(z)^{-1}$ is analytic for $z \in {\mathcal U}_{s,r}$ then Corollary \ref{fundeqnunique} shows that the coefficients are uniquely defined and Theorem \ref{fundsol} shows that
\begin{eqnarray*}
R(z) & = & \sum_{p=1}^{\infty} (-1)^{p-1}(R_{-1}A_0)^{p-1}R_{-1}z^{-p} + \sum_{n=0}^{\infty}
(-1)^n(R_0A_1)^nR_0z^n \\
& = & (I z + R_{-1}A_0)^{-1}R_{-1} + (I + R_0A_1z)^{-1}R_0.
\end{eqnarray*}

Conversely suppose there exist operators $R_{-1}, R_0 \in {\mathcal L}(K,H)$ such that conditions (i)--(iii) in the statement of the corollary are satisfied.  Define $R:{\mathcal U}_{s,r} \rightarrow {\mathcal L}(K,H)$ as above.  Since $R_{-1}A_0R_0 = 0$ we have
\begin{eqnarray*}
\lefteqn{ (I z + R_{-1}A_0)^{-1}R_{-1}A(z) } \\
& = &  \sum_{p=1}^{\infty} (-1)^{p-1}(R_{-1}A_0)^{p-1}R_{-1}z^{-p} (A_0 + A_1z) \\
& = &  \sum_{p=1}^{\infty} (-1)^{p-1}(R_{-1}A_0)^pz^{-p} + \sum_{p=1}^{\infty} (-1)^{p-1}(R_{-1}A_0)^{p-1}(I -
R_0A_0)z^{-p+1} \\
& = &  \sum_{p=1}^{\infty} (-1)^{p-1}(R_{-1}A_0)^pz^{-p} + \sum_{p=1}^{\infty}
(-1)^{p-1}(R_{-1}A_0)^{p-1}z^{-p+1} - R_0A_0 \\
& = & P
\end{eqnarray*}
for $|z| > s$.  Since $R_0A_1R_{-1} = 0$ it follows that
\begin{eqnarray*}
\lefteqn{ \hspace{-0.5cm} (I + R_0A_1z)^{-1}R_0A(z) } \\
& = & \sum_{n=0}^{\infty} (-1)^n(R_0A_1)^nR_0z^n (A_0 + A_1z) \\
& = & \sum_{n=0}^{\infty} (-1)^n(R_0A_1)^n(I - R_{-1}A_1) z^n +  \sum_{n=0}^{\infty} (-1)^n(R_0A_1)^{n+1} z^{n+1}
\\
& = & \sum_{n=0}^{\infty} (-1)^n(R_0A_1)^nz^n - R_{-1}A_1 +  \sum_{n=0}^{\infty} (-1)^n(R_0A_1)^{n+1} z^{n+1} \\
& = & I - P
\end{eqnarray*}
for $|z| < r$.  Thus $R(z)A(z) = I$ for $z \in {\mathcal U}_{s,r}$.  A similar argument shows $A(z)R(z) = I$
for $z \in {\mathcal U}_{s,r}$.  Hence, from Corollary \ref{fundeqnunique}, the coefficients are uniquely defined.  From (\ref{cf}),
\begin{eqnarray*}
\lefteqn{ \hspace{-0.3cm} {\mathcal R}_{\lambda} - {\mathcal R}_{\mu} }  \\
& = & (A_0 \lambda - A_1)^{-1}A_0 - (A_0 \mu - A_1)^{-1}A_0 \\
& = & - (I - R_{-1}A_0\lambda)^{-1}R_{-1}A_0 + \lambda^{-1}(I - R_0A_1\lambda^{-1})^{-1}R_0A_0 \\
& & \hspace{0.5cm} + (I - R_{-1}A_0\mu)^{-1}R_{-1}A_0 - \mu^{-1}(I - R_0A_1\mu^{-1})^{-1}R_0A_0 \\
& = & (\mu - \lambda) (R_{-1}A_0)^2 \left[ I - R_{-1}A_0(\lambda + \mu) + (R_{-1}A_0)^2(\lambda^2 + \lambda \mu +
\mu^2) - \cdots \right] \\
& & \hspace{0.5cm} + (\lambda^{-1} - \mu^{-1}) \left[ I - R_0A_1(\lambda^{-1} + \mu^{-1}) \right. \\
& & \hspace{1cm} \left. + (R_0A_1)^2(\lambda^{-2} + \lambda^{-1}\mu^{-1} + \mu^{-2}) - \cdots \right]R_0A_0 \\
& = & (\mu - \lambda) \left\{ \rule{0cm}{0.4cm} \left[ I - R_1A_0 \lambda + \cdots \right] R_{-1}A_0 \cdot
\left[ I - R_1A_0 \mu + \cdots \right] R_{-1}A_0 \right. \\
& & \hspace{0.5cm} \left. + \lambda^{-1} \left[ I - \lambda^{-1}R_0A_1 + \cdots \right]R_0A_0 \cdot \mu^{-1} \left[
I - \mu^{-1}R_0A_1 + \cdots \right]R_0A_0 \rule{0cm}{0.4cm} \right\} ,
\end{eqnarray*}
where we have used the fact that $(R_0A_0)^2 = R_0A_0$ and $(R_1A_0)R_0A_0 = - (R_0A_1)R_0A_0 = - R_0(A_1R_0)A_0 = R_0(A_0R_1)A_0$.  Therefore
\begin{eqnarray*}
{\mathcal R}_{\lambda} - {\mathcal R}_{\mu} & = & (\mu - \lambda) \left\{ \rule{0cm}{0.4cm} (I -
R_{-1}A_0\lambda)^{-1}R_{-1}A_0 \cdot (I - R_{-1}A_0\mu)^{-1}R_{-1}A_0 \right. \\
& & \hspace{0.5cm} \left. + \lambda^{-1}(I - R_0A_1\lambda^{-1})^{-1}R_0A_0 \cdot  \mu^{-1}(I -
R_0A_1\mu^{-1})^{-1}R_0A_0 \rule{0cm}{0.4cm} \right\} \\
& = & (\mu - \lambda) {\mathcal R}_{\lambda} {\mathcal R}_{\mu}
\end{eqnarray*}
where the final step depends on the identities $R_{-1}A_0R_0 = 0$ and $R_0A_0R_{-1} = 0$.  A similar argument
establishes the other resolvent equation.  $\hfill \Box$

\begin{definition}
\label{bs}
If the coefficients $\{ R_{-1},R_0 \} \in {\mathcal L}(K,H)$ satisfy $(i)$--$(iii)$ in Corollary $\ref{basicsol}$ then we say $\{R_{-1},R_0\}$ is the basic solution to $(\ref{lfe})$ and $(\ref{rfe})$ on ${\mathcal U}_{s,r}$.
\end{definition}

\begin{remark}
\label{formulation}
Although one could argue retrospectively that existence of an inverse operator implies that the spaces $H$ and $K$ must be isomorphic there are many problems where they are best regarded as distinct.  See for instance {\em \cite{how4}} where examples are given involving perturbed bounded linear differential operators mapping a particular Sobolev space $H$ onto a space $K$ of square integrable functions.
\end{remark}

\section{Isolated singularities}
\label{iso}

Choose $j \in {\mathbb N}$.  We note that $R_{-j} \neq 0 \iff (R_{-1}A_0)^{j-1} \neq 0 \iff (A_0R_{-1})^{j-1} \neq 0$.  From the properties in Theorem \ref{fundsol} we see, on the one hand, that $R_{-j} = (-1)^{j-1}(R_{-1}A_0)^{j-1}R_{-1} = (-1)^{j-1}R_{-1}(A_0R_{-1})^{j-1}$ and, on the other hand, that $(-1)^{j-1}(R_{-1}A_0)^{j-1} = R_{-j}A_1$ and $(-1)^{j-1}(A_0R_{-1})^{j-1} = A_1R_{-j}$.  Hence we can use the properties of either $R_{-1}A_0 \in {\mathcal L}(H)$ or $A_0R_{-1} \in {\mathcal L}(K)$ to describe an isolated singularity of $R(z)$ at $z=0$.  If $(R_{-1}A_0)^j \neq 0$ for $j < p$ and $(R_{-1}A_0)^j = 0$ for $j \geq p$, then there is a pole of order $p$.  Clearly $s = \lim_{j \rightarrow \infty}  \|(R_{-1}A_0)^j\|^{1/j} = 0$ in this case.  Alternatively, if $s = \lim_{j \rightarrow \infty}  \|(R_{-1}A_0)^j\|^{1/j} = 0$ but $(R_{-1}A_0)^j \neq 0$ for all $j \in {\mathbb N}$, then there is an essential singularity.  Our results extend classical resolvent theory for isolated singularities \cite[pp. 229--231 and pp. 282--286]{yos1} to linear pencils.

In the next section we show how Gaussian elimination can be used to obtain a direct solution to the fundamental equations for the resolvent of a matrix operator near an isolated singularity.  A general solution procedure for the fundamental equations remains an open question.  For matrix operators some authors  \cite{avr1,wil1} have used Jordan chains to find the Laurent series.  We discuss Jordan chains briefly in the next two subsections but a full discussion is beyond the scope of this paper.  We refer readers to \cite{wil1} for a more extensive discussion.  

\subsection{Jordan chains for a finite-order pole}

Following Gohberg and Rodman \cite{goh2} we say that $\bfx_0,\ldots,\bfx_{p-1}$ is a Jordan chain of length $p$ for $A(z)$ if $A_0 \bfx_0 = \bfzero$ and $A_1\bfx_{j-1} + A_0\bfx_j = \bfzero$ for $j=1,\ldots,p-1$.  If $R(z)$ has a pole of order $p$ at $z=0$ then $(R_{-1}A_0)^{p-1} R_{-1} \neq 0$ and hence there must be some $\bphi \in K$ with $(R_{-1}A_0)^{p-1} R_{-1} \bphi \neq \bfzero$.   If we define $\bfx_j = (-1)^{p-j-1}(R_{-1}A_0)^{p-j-1}R_{-1}\bphi$ for all $j=0,\ldots,p-1$ then
$$
A_0 \bfx_0 = (-1)^{p-1}(A_0R_{-1})^p \bphi = \bfzero
$$
and
\begin{eqnarray*}
A_1\bfx_{j-1} + A_0\bfx_j & = & A_1(-1)^{p-j}(R_{-1}A_0)^{p-j} R_{-1} \bphi + A_0(-1)^{p-j-1}(R_{-1}A_0)^{p-j-1}R_{-1} \bphi \\
& = & (-1)^{p-j}[I - A_0R_0](A_0R_{-1})^{p-j} \bphi + (-1)^{p-j-1}(A_0R_{-1})^{p-j} \bphi \\
& = & 0
\end{eqnarray*}
for each $j=1,\ldots,p-1$ since $R_0A_0R_{-1} = 0$.  Hence $\bfx_0,\ldots,\bfx_{p-1}$ forms a Jordan chain of length $p$ for $A(z)$.  For matrix operators $A_0, A_1 \in {\mathbb C}^{n \times n}$ where $A_0$ is singular then $R(z) = \mbox{adj}\, A(z)/ \det A(z)$ and hence $R(z)$ has a pole of order $p \leq n$ at $z = 0$.   

\subsection{Jordan chains for an isolated essential singularity}  

First we show that an isolated essential singularity is possible.  Choose $\delta \in (0,1)$ and define $W = [w_{i,j}] \in {\mathcal L}(\ell_2)$ by setting $w_{i,i+1} = \delta^{2^{i-1}}$ for all $i \in {\mathbb N}$ with $w_{i,j} = 0$ otherwise.  Therefore $W^n = [w_{i,j}^{(n)}]$ is given by $w_{i,i+n}^{(n)} = \delta^{\kappa(i,n)}$ where $\kappa(i,n) = 2^{i-1}(2^n - 1)$ for all $i \in {\mathbb N}$ with $w_{i,j}^{(n)} = 0$ otherwise.  Since $W^n \bfe_{n+1} = \delta^{2^n-1} \bfe_1$ it follows that $\|W^n\| = \delta^{2^n-1}$ and hence that $\|W^n\|^{1/n} \rightarrow 0$ as $n \rightarrow \infty$.  Note that $W^n \neq 0$ for all $n \in {\mathbb N}$.  If $A(z) = Iz - W$ then
\begin{equation}
\label{isex}
(I z - W)^{-1} = \frac{1}{z} \left[ I + \frac{W}{z} + \left( \frac{W}{z} \right)^2 + \cdots \right]
\end{equation}
for $z \neq 0$.  Thus $R(z) = (I z - W)^{-1}$ has an essential singularity at $z = 0$.  We will not discuss Jordan chains in general but we make the following observations.  In this example we have $A_0 = - W$, $A_1 = I$, $R_{-1} = I$ and $R_0 = 0$.  Hence $R_{-1}A_0 = - W$.  Since $W \bfe_{k+1} = \delta^{2^{k-1}} \bfe_k$ for each $k \in {\mathbb N}$ and since $W \bfe_1 = \bfzero$ it follows that if $\bphi = \bfe_{k+1}$ then $W^{k+1} \bphi = \bfzero$ and $\bphi, W \bphi,\ldots,W^k \bphi$ forms a Jordan chain of length $k+1$.  Thus we can find a Jordan chain of any finite length for this example.  Note also that if $\bphi = \sum_{k \in {\mathbb N}} \bfe_k/k$ then $W^n \bphi \neq \bfzero$ for all $n \in {\mathbb N}$. 

\section{Solving the fundamental equations}
\label{sfe}

In the case of a pole of order $p$ it has been shown \cite{alb1} that a solution for $\{R_{-p},\ldots,R_0\}$ which includes the basic solution $\{R_{-1},R_0\}$ can be determined by solving a finite subset of the fundamental equations---the {\em so-called} determining equations.  The order of the pole and the solution can be obtained by applying block row reductions to the augmented operator
\begin{equation}
\label{determining}
\left[ \begin{array}{cccc|cccc}
A_0 & 0 & 0 & \cdots & I & 0 & 0 & \cdots \\
A_1 & A_0 & 0 & \cdots & 0 & I & 0 & \cdots \\
0 & A_1 & A_0 & \cdots & 0 & 0 & I  & \cdots \\
\vdots & \vdots & \vdots &  & \vdots & \vdots & \vdots &  \end{array} \right].
\end{equation}
The first column on the right-hand side gives a solution if $A(z)$ is analytic at $z=0$ but is inconsistent otherwise; and for $p \geq 0$ the $(p+1)^{\mbox{\scriptsize th}}$ column on the right-hand side gives a solution if $A(z)$ has a pole of order $p$ but is inconsistent if there is a higher-order pole.  For a matrix operator $A(z) \in {\mathbb C}^{m \times m}$, the solution procedure for the determining equations (\ref{determining}) is simply a standard Gaussian elimination.  Note that the special structure of the augmented coefficient matrix means the row reduction is recursive and hence can be conveniently implemented for numerical calculations.  We consider a simple example.

\begin{example}
\label{p2ex}

Let
$$
A_0 = \left[ \begin{array}{ccc}
1 & 0 & 1 \\
1 & 0 & 0 \\
1 & 0 & 0 \end{array} \right] \quad \mbox{and} \quad A_1 = \left[ \begin{array}{ccc}
1 & 0 & -1 \\
0 & 1 & 0 \\
0 & 1 & 1 \end{array} \right].
$$
Gaussian elimination on the system $(\ref{determining})$ shows a pole of order $2$ with solution
$$
R_{-2} = \left[\begin{array}{ccc}
0 & 0 & 0 \\
0 & -1 & 1 \\
0 & 0 & 0 \end{array} \right], \quad R_{-1} = \left[\begin{array}{ccc}
0 & 1 & -1 \\
-1 & 3 & -2 \\
0 & -1 & 1 \end{array} \right]\quad R_{0} = \left[\begin{array}{ccc}
1 & -2 & 2 \\
1 & -2 & 2 \\
0 & 0 & 0 \end{array} \right].
$$
\end{example}
For operators on Hilbert space with a finite-order pole the same finite system of equations is sufficient but the block row operations must be considered more generally---even if an appropriate orthonormal basis of eigenvectors can be conveniently used for an infinite matrix formulation.  In general, for Banach space and for isolated essential singularities there is no known guaranteed standard solution procedure.

\section{Unbounded operators}

If the unperturbed operator $A_0:{\mathcal D} \subset H \rightarrow K$ is densely defined and closed but unbounded we may use a Sobolev space argument to modify the topology and convert $A_0$ into a bounded operator.  The operator $A_0$ is closed if and only if $\{\bfx_n\}_{n \in {\mathbb N}} \in {\mathcal D}$ with $\|\bfx_n - \bfx\| \rightarrow 0$ and $\|A_0 \bfx_n - \bfy\| \rightarrow 0$ as $n \rightarrow \infty$ for $\bfx \in H$ and $\bfy \in K$ implies $\bfx \in {\mathcal D}$ and $A_0 \bfx = \bfy$.  In such cases we may define a new Banach space
$$
H_E = \{ \bfx \in {\mathcal D} \mid \|\bfx\| + \|A_0\bfx\| < \infty\}
$$
with $\| \bfx\|_E = \|\bfx\| + \|A_0\bfx\|$.  The operator $A_{0,E}:H_E \rightarrow K$ defined by $A_{0,E}\bfx = A_0\bfx$ for $\bfx \in H_E$ is bounded with
$$
\|A_{0,E}\| = \sup_{\sbfx \in H_E} \frac{\|A_{0,E} \bfx\|}{\| \bfx\|_E} = \sup_{\sbfx \in H_E} \frac{\|A_0 \bfx\|}{\|\bfx\| + \|A_0 \bfx\|}  \leq 1.
$$
The differential operator in one or more dimensions provides an important example of a linear operator which is densely defined and closed but nevertheless unbounded.  It is standard practice in modern analysis to treat this operator as a bounded operator on an appropriate Sobolev space.  Consider the following example.

\begin{example}[Gradient approximation]
\label{grad}

Let ${\mathcal U} = [0,1]^3$ denote the unit cube in ${\mathbb R}^3$ with boundary $\partial {\mathcal U}$.  Suppose an observed square integrable function $\bfg:{\mathcal U} \rightarrow {\mathbb R}^3$ is given.  We wish to find the potential function $u:{\mathcal U} \rightarrow {\mathbb R}$ with $u(\bfx) = 0$ when $\bfx \in \partial {\mathcal U}$, which minimises the total residual error
$$
\int_{\mathcal U} \| \bnab u(\bfx) - \bfg(\bfx) \|^2\ d\bfx.
$$
Define the Hilbert space $H$ of functions $u:{\mathcal U} \rightarrow {\mathbb R}$ such that
$$
\int_{\mathcal U} \left[ \rule{0cm}{0.4cm}\ |u(\bfx)|^2 + \|\bnab u(\bfx) \|^2\ \right]\ d\bfx < \infty
$$
with $u(\bfx) = 0$ when $\bfx \in \partial {\mathcal U}$ and inner product
$$
\langle u,v \rangle_H = \int_{\mathcal U} \left[ \rule{0cm}{0.4cm}\ u(\bfx)v(\bfx) + \langle \bnab u(\bfx),\bnab v(\bfx) \rangle\ \right] d\bfx
$$
for each $u,v \in H$.  Let $K$ be the Hilbert space of square integrable functions $\bfg:{\mathcal U} \mapsto {\mathbb R}^3$.  The gradient operator $\bnab:H \rightarrow K$ is a bounded linear map.  For more information and a solution see {\em \cite{how2}}.  The solution involves a bounded self-adjoint differential operator.

\end{example}

Other examples involving perturbed operators can be found in \cite{how4}.

\section{Global structure of the resolvent}
\label{gsr}

In order to calculate the resolvent it is convenient to represent the spaces $H$ and $K$ as direct sums of projected subspaces.  We can use the natural corresponding projections $P = R_{-1}A_1 \in {\mathcal L}(H)$ and $Q = A_1R_{-1} \in {\mathcal L}(K)$ to write $H = M \oplus M^c$ and $K = N \oplus N^c$ where $M = P(H)$, $M^c = (I-P)(H)$, $N = Q(K)$ and $N^c = (I-Q)(K)$.  Thus we can write the local spectral decomposition in the form
\begin{equation}
\label{hklocaldecomp}
R(z) = R_{[0]}(z) + R_{[0]}^c(z)
\end{equation}
where $R_{[0]}: N \rightarrow M$ and $R_{[0]}^c:N^c \rightarrow M^c$ are defined respectively by $R_{[0]}(z) = PR(z)Q = (Iz + R_{-1}A_0)^{-1}R_{-1}$ for $|z| > s$ which represents the bounded component of the spectral set (the singularities with $|z| \leq s$) and $R_{[0]}^c(z) = (I-P)R(z)(I-Q) = (I +R_0A_1z)^{-1}R_0$ for $|z| < r$ which represents the unbounded component (the singularities with $|z| \geq r$).

Suppose $R(z)$ has isolated singularities at $z=z_k$ for each $k =1,2,\ldots,m$ and is analytic elsewhere.   Write $A(z) = A_{k,0} + A_{k,1}(z-z_k)$ where $A_{k,0} = A_0 + A_1z_k$ and $A_{k,1} = A_1$.   Let $\{ R_{k,-1}, R_{k,\,0} \}$ denote a basic solution to the fundamental equations on $z_k +\ {\mathcal U}_{0,\,r_k}$ where $r_k = \min_{j \neq k} |z_j - z_k|$.  Let $P_k = R_{k,-1}A_1$ and $Q_k = A_1R_{k,-1}$ be the corresponding projections and write $R(z) = R_{[k]}(z) + R_{[k]}^c(z)$ where
\begin{equation}
\label{ksing}
R_{[k]}(z) = P_kR(z)Q_k \quad (z \neq z_k)
\end{equation}
is singular at $z=z_k$ and where
\begin{equation}
\label{kreg}
R_{[k]}^c(z) = (I-P_k)R(z)(I-Q_k) \quad (z \neq z_1,\ldots,z_{k-1},z_{k+1},\ldots,z_m)
\end{equation}
is regular at $z=z_k$.  We now show that $P_kP_{\ell} = 0$ and $Q_kQ_{\ell} = 0$ for $k \neq \ell$.  Note from the proof of Corollary \ref{basicsol} that $R_{[k]}(z)A(z) = P_k$ for $z \neq z_k$ and $A(z)R_{[\ell]}(z) = Q_{\ell}$ for $z \neq z_{\ell}$.  Firstly
\begin{equation}
\label{pkrell}
R_{[k]}(z)A(z)R_{[\ell]}(z) = P_k R_{[\ell]}(z) = \sum_{p=0}^{\infty} \frac{P_k {R_{[\ell]}}^{(p)}(z_k)}{p!} (z - z_k)^p
\end{equation}
for $|z - z_k| < |z_{\ell} - z_k|$ since $R_{[\ell]}(z)$ is analytic at $z = z_k$ and secondly
\begin{equation}
\label{rkqell}
R_{[k]}(z)A(z)R_{[\ell]}(z) = R_{[k]}(z)Q_{\ell} = \sum_{q=1}^{\infty} \frac{R_{k,-q}Q_{\ell}}{(z - z_k)^q}
\end{equation}
for $z \neq z_k$.  Since (\ref{pkrell}) and (\ref{rkqell}) are equal for $0 < |z - z_k| < |z_{\ell} - z_k|$ it follows by equating coefficients that $P_k {R_{[\ell]}}^{(p)}(z_k) = 0$ for all $p \in \{0\} \cup {\mathbb N}$ and $R_{k,-q}Q_{\ell} = 0$ for all $q \in {\mathbb N}$.  Therefore $R_{[k]}(z)A(z)R_{[\ell]}(z) = 0$ for $0 < |z - z_k| < |z_{\ell} - z_k|$.  In particular $P_kR_{[\ell]}(z_k) = 0 \Rightarrow R_{k,-1}A_1R_{\ell,-1} \left[ I(z_k - z_{\ell}) + A_{\ell,0}R_{\ell,-1} \right]^{-1} = 0 \Rightarrow R_{k,-1}A_1R_{\ell,-1} = 0$ and so $P_kP_{\ell} = 0$ and $Q_kQ_{\ell} = 0$ for $k \neq \ell$.  Hence there also exist corresponding complementary projections $P_{\infty} = \prod_{k=1}^m (I - P_k) = I - \sum_{k=1}^m P_k \in {\mathcal L}(H)$ and $Q_{\infty} = \prod_{k=1}^m (I - Q_k) = I - \sum_{k=1}^m Q_k \in {\mathcal L}(K)$.  Thus we may write
\begin{equation}
\label{hkdecomp}
H = M_1 \oplus \cdots \oplus M_m \oplus M_{\infty} \quad \mbox{and} \quad K = N_1 \oplus \cdots \oplus N_m
\oplus N_{\infty}
\end{equation}
where $M_k = P_k(H)$, $N_k = Q_k(K)$, $M_{\infty} = P_{\infty}(H) = (I - \sum_{k=1}^m P_k)(H)$ and $N_{\infty} = Q_{\infty}(K) = (I - \sum_{k=1}^m Q_k)(K)$ and thereby establish the full spectral decomposition
\begin{equation}
\label{gs}
R(z) = \sum_{k=1}^m R_{[k]}(z) + R_{[\infty]}(z)
\end{equation}
where each term $R_{[k]}(z) = P_kR(z)Q_k \in {\mathcal L}(N_k,M_k)$ is analytic for $z \neq z_k$ and the remainder $R_{[\infty]}(z) = P_{\infty}R(z)Q_{\infty} \in {\mathcal L}(N_{\infty},M_{\infty})$ is entire.  We illustrate our result with an elementary example.

\begin{example}
\label{globalex}

Let $H = K = {\mathbb C}^3$ and define $A_0,A_1 \in {\mathcal L}(H,K)$ by
$$
A_0 = \left[ \begin{array}{ccc}
1 & 1 & 1 \\
1 & 2 & 1 \\
2 & 1 & 2 \end{array} \right] \quad \mbox{and} \quad A_1 = \left[ \begin{array}{ccc}
1 & 0 & 0 \\
0 & 1 & 0 \\
1 & 0 & 1 \end{array} \right].
$$
The resolvent of $A(z) = A_0 + A_1z$ is singular at $z_1=0$, $z_2 = -1$ and $z_3= -3$.  Hence it can be shown that
$$
R(z) = \frac{R_{1,-1}}{z} + \frac{R_{2,-1}}{z+1} + \frac{R_{3,-1}}{z+3}
$$
for $z \neq 0, -1, -3$ where
$$
R_{1,-1} = \left[ \begin{array}{ccc}
\rule{0cm}{0.45cm} 1 & -\, \frac{1}{3} & -\, \frac{1}{3} \\
\rule{0cm}{0.45cm} 0 & 0 & 0 \\
\rule{0cm}{0.45cm} -1 & \frac{1}{3} & \frac{1}{3} \end{array} \right], \quad R_{2,-1} = \left[ \begin{array}{ccc}
\rule{0cm}{0.45cm} 0 & 0 & 0 \\
\rule{0cm}{0.45cm} 0 & \frac{1}{2} & - \frac{1}{2} \\
\rule{0cm}{0.45cm} 0 & - \frac{1}{2} & \frac{1}{2} \end{array} \right], \quad R_{3,-1} =  \left[ \begin{array}{ccc}
\rule{0cm}{0.45cm} 0 & \frac{1}{3} & \frac{1}{3} \\
\rule{0cm}{0.45cm} 0 & \frac{1}{2} & \frac{1}{2} \\
\rule{0cm}{0.45cm} 0 & \frac{1}{6} & \frac{1}{6} \end{array} \right].
$$
To obtain Laurent series in each of the relevant annular regions centred at $z=0$ we simply use the well-known expansions
$$
1/(z+a) = (1/a)[ 1 - (z/a) +(z/a)^2 - (z/a)^3 + \cdots \,]
$$
for $|z| < |a|$ and
$$
1/(z+a) = (1/z)[ 1 - (a/z) + (a/z)^2 - (a/z)^3 + \cdots \,]
$$
for $|z| > |a|$.  Simple calculations give the following results.  In the region $0 < |z| < 1$ we have
$$
R_{-1} = \left[ \begin{array}{ccc}
\rule{0cm}{0.45cm} 1 & -\frac{1}{3} & - \frac{1}{3} \\
\rule{0cm}{0.45cm} 0 & 0 & 0 \\
\rule{0cm}{0.45cm} -1 & \frac{1}{3} & \frac{1}{3} \end{array} \right] \quad \mbox{and} \quad R_0 = \left[ \begin{array}{ccc}
\rule{0cm}{0.45cm} 0 & \frac{1}{9} & \frac{1}{9} \\
\rule{0cm}{0.45cm} 0 & \frac{2}{3} & -\frac{1}{3} \\
\rule{0cm}{0.45cm} 0 & -\frac{4}{9} & \frac{5}{9} \end{array} \right]
$$ 
from which it follows that $\|R_{-1}A_0\| = 0$ and $\|R_0A_1\| = 1$ as expected.  In the region $1 < |z| < 3$ we obtain
$$
R_{-1} = \left[ \begin{array}{ccc}
\rule{0cm}{0.45cm} 1 & -\frac{1}{3} & - \frac{1}{3} \\
\rule{0cm}{0.45cm} 0 & \frac{1}{2} & - \frac{1}{2} \\
\rule{0cm}{0.45cm} -1 & -\frac{1}{6} & \frac{5}{6} \end{array} \right] \quad \mbox{and} \quad R_0 = \left[ \begin{array}{ccc}
\rule{0cm}{0.45cm} 0 & \frac{1}{9} & \frac{1}{9} \\
\rule{0cm}{0.45cm} 0 & \frac{1}{6} & \frac{1}{6} \\
\rule{0cm}{0.45cm} 0 & \frac{1}{18} & \frac{1}{18} \end{array} \right]
$$
which gives $\| (R_{-1}A_0)^j \| = 1$ for all $j \in {\mathbb N}$ and $\| R_0A_1 \| = 1/3$ in line with the theory.  For the final region $3 < |z| < \infty$ the calculations show that
$$
R_{-1} = \left[ \begin{array}{ccc}
\rule{0cm}{0.45cm} 1 & 0 & 0 \\
\rule{0cm}{0.45cm} 0 & 1 & 0  \\
\rule{0cm}{0.45cm} -1 & 0 & 1 \end{array} \right] \quad \mbox{and} \quad R_0 = \left[ \begin{array}{ccc}
\rule{0cm}{0.45cm} 0 & 0 & 0 \\
\rule{0cm}{0.45cm} 0 & 0 & 0 \\
\rule{0cm}{0.45cm} 0 & 0 & 0 \end{array} \right]
$$
from which it follows that $\|(R_{-1}A_0)^j\| = 3^j$ as required.  It is easy to check that the fundamental equations are satisfied on each region.  

\end{example}

\section{Polynomial pencils}
\label{pp}

Finally, let $A_0, \ldots, A_{n} \in {\mathcal L}(H,K)$, where $A_0$ is singular.  Define a polynomial pencil $A: {\mathbb C} \rightarrow {\mathcal L}(H,K)$ by the formula $A(z) = \sum_{i=0}^{n} A_i z^i$ and consider a resolvent $R:{\mathcal U}_{s,r} \rightarrow {\mathcal L}(K,H)$ for $A$ defined on an annulus ${\mathcal U}_{s,r}$ for some $0 \leq s < r$ by a Laurent series $R(z) = \sum_{j \in {\mathbb Z}} R_jz^j$ where $R_j \in {\mathcal L}(K,H)$ for each $j \in {\mathbb Z}$.  Consider also a corresponding augmented linear operator pencil ${\mathcal A}:{\mathbb C} \rightarrow {\mathcal L}(H^n,K^n)$ defined by the formula ${\mathcal A}(z) = {\mathcal A}_0 + {\mathcal A}_1 z$ where
\begin{equation}
\label{augop}
{\mathcal A}_0 = \left[ \begin{array}{cccc}
A_0 & 0 & \cdots & 0 \\
A_1 & A_0 & \cdots & 0 \\
\vdots & \vdots & \ddots & \vdots \\
A_{n-1} & A_{n-2} & \cdots & A_0 \end{array} \right]  \quad \mbox{and} \quad {\mathcal A}_1 = \left[ \begin{array}{cccc}
A_{n} & A_{n-1} & \cdots & A_1 \\
0 & A_{n} & \cdots & A_2 \\
\vdots & \vdots & \ddots & \vdots \\
0 & 0 & \cdots & A_{n} \end{array} \right].
\end{equation}
Define an augmented operator ${\mathcal R}:{\mathcal U}_{s,r} \rightarrow {\mathcal L}(K^n,H^n)$ by the Laurent series ${\mathcal R}(z) = \sum_{j \in {\mathbb Z}} {\mathcal R}_jz^j$ where ${\mathcal R}_j \in {\mathcal L}(K^n,H^n)$ is defined by the formula
\begin{equation}
\label{augr}
{\mathcal R}_j = \left[ \begin{array}{cccc}
R_{nj} & R_{nj-1} & \cdots & R_{nj-n+1} \\
R_{nj+1} & R_{nj} & \cdots & R_{nj-n+2} \\
\vdots & \vdots & \ddots & \vdots \\
R_{nj+n-1} & R_{nj+n-2} & \cdots & R_{nj} \end{array} \right]
\end{equation}
for each $j \in {\mathbb Z}$.   It is a matter of elementary algebra to show that the left fundamental equations
\begin{equation}
\label{plfe}
\sum_{k=0}^n R_{j-n+k}A_{n-k} = \left\{ \begin{array}{ll}
I & \mbox{if}\ j = 0 \\
0 & \mbox{if}\ j \neq 0 \end{array} \right.
\end{equation}
and right fundamental equations
\begin{equation}
\label{prfe}
\sum_{k=0}^n A_{n-k}R_{j-n+k} = \left\{ \begin{array}{ll}
I & \mbox{if}\ j = 0 \\
0 & \mbox{if}\ j \neq 0 \end{array} \right.
\end{equation}
for the polynomial pencil are equivalent to the left fundamental equations
\begin{equation}
\label{auglfe}
{\mathcal R}_{j-1}{\mathcal A}_1 + {\mathcal R}_j{\mathcal A}_0 = \left \{ \begin{array}{ll}
{\mathcal I} & \mbox{if}\ j=0 \\
0 & \mbox{if}\ j \neq 0 \end{array} \right.
\end{equation}
and right fundamental equations
\begin{equation}
\label{augrfe}
{\mathcal A}_1{\mathcal R}_{j-1} + {\mathcal A}_0{\mathcal R}_j = \left \{ \begin{array}{ll}
{\mathcal I} & \mbox{if}\ j=0 \\
0 & \mbox{if}\ j \neq 0 \end{array} \right.
\end{equation}
for the augmented linear pencil.  In (\ref{auglfe}) and (\ref{augrfe}) we have used the symbol ${\mathcal I}$ to denote the identity operator on both $H^n$ and $K^n$.  It is a simple matter to find corresponding geometric bounds on the coefficients of the respective resolvents.   Thus the resolvent $R:{\mathcal U}_{s,r} \rightarrow {\mathcal L}(K,H)$ for the polynomial pencil $A$ is well defined if and only if the corresponding augmented resolvent ${\mathcal R}:{\mathcal U}_{s,r} \rightarrow {\mathcal L}(K^n,H^n)$ for the augmented linear pencil ${\mathcal A}$ is well defined.  Hence our results for linear pencils can also be applied to polynomial pencils. 

We will not write down a complete list of analogous results for polynomial pencils but rather be content to illustrate our assertions by considering the concept of a basic solution.  The suitably bounded pair $\{{\mathcal R}_{-1}, {\mathcal R}_0 \}$ provides a basic solution to the fundamental equations for the augmented linear pencil in that the complete solution to the fundamental equations (\ref{auglfe}) and (\ref{augrfe}) is given by the set $\{ {\mathcal R}_j \}_{j \in {\mathbb Z}}$ where
\begin{equation}
\label{augcoeff}
{\mathcal R}_j = \left \{ \begin{array}{ll}
(-1)^j ({\mathcal R}_0{\mathcal A}_1)^j {\mathcal R}_0 & \mbox{for}\ j  \geq 0 \\
& \\
(-1)^{-j-1} ({\mathcal R}_{-1}{\mathcal A}_0)^{-j-1}{\mathcal R}_{-1} & \mbox{for}\ j < 0 \end{array} \right.
\end{equation}
for each $j \in {\mathbb Z}$.  Thus the set of coefficients $\{R_j\}_{j \in {\mathbb Z}}$ for the resolvent $R$ of the polynomial pencil is completely defined by $\{ {\mathcal R}_{-1}, {\mathcal R}_0 \}$ and hence also by $\{R_{-2n+1}, R_{-2n+2},\ldots,R_{n-1}\}$.  Thus we have the following definition.

\begin{definition}
\label{augbs}
We say that $\{R_{-2n+1}, R_{-2n+2}, \ldots, R_{n-1}\} \in {\mathcal L}(K,H)$ is a basic solution\footnote{Note that the basic solution is wrongly listed as $\{R_{-n}, R_{-n+1},\ldots,R_n\}$ in the published article.} to $(\ref{plfe})$ and $(\ref{prfe})$ on ${\mathcal U}_{s,r}$ if and only if $\{{\mathcal R}_{-1},{\mathcal R}_0 \} \in {\mathcal L}(K^n,H^n)$ is a basic solution to $(\ref{auglfe})$ and $(\ref{augrfe})$ on ${\mathcal U}_{s,r}$.
\end{definition}

The correspondence between the resolvents for the polynomial pencil and the augmented linear pencil can be used in an analogous way to obtain results for polynomial pencils corresponding to Theorems \ref{fundeqn} and \ref{fundsol} and to Corollary \ref{basicsol}.

\section{Analytic pencils}
\label{ap}

For analytic pencils where the coefficients satisfy a finite recursion it has been shown \cite{how1} that inversion of the analytic pencil is equivalent to inversion of a polynomial pencil.  However we also note more generally that analytic pencils can be inverted if they contain a polynomial pencil that can be inverted.  Let $A(z) = \sum_{j=0}^{\infty} A_jz^j$ where $A_j \in {\mathcal L}(H,K)$ be an analytic pencil where $A_0$ is singular.  Suppose the resolvent $R(z) = A(z)^{-1}$ is defined for $z \in {\mathcal U}_{s,r}$.  From the Banach inverse theorem it follows that for each $z \in {\mathcal U}_{s,r}$ there is some $\alpha(z) > 0$ such that
\begin{equation}
\label{ap1}
\alpha(z) = \inf_{\sbfx \in H} \|A(z) \bfx \| / \|\bfx\|
\end{equation}
 for all $\bfx \in H$.  Choose $\delta, \epsilon \in {\mathbb R}$ such that $s < s + \epsilon < r - \delta < r$.  Since $\alpha(z) >0$ is continuous on the compact set $F_{s+\epsilon,r-\delta} = \{ z \in {\mathbb C} \mid s+\epsilon \leq |z| \leq r - \delta\}$ it follows that there is some $\alpha > 0$ such that $\alpha(z) \geq \alpha$ for all $z \in F_{s+\epsilon,r-\delta}$.  Write
 $$
 A(z) = \sum_{j=0}^m A_jz^j + \sum_{j=m+1}^{\infty} A_jz^j = A_m(z) + \rho_m(z)
 $$
 where $A_m(z)$ is the partial sum and $\rho_m(z)$ is the remainder.  Since $\|A(z) - A_m(z)\| \rightarrow 0$ as $m \rightarrow \infty$ uniformly on $F_{s+\epsilon,r-\delta}$ it follows that we can find $m$ such that $\|A_m(z) \bfx \| \geq (\alpha/2) \|\bfx\|$ for all $z \in F_{s+\epsilon,r-\delta}$ and all $\bfx \in H$.   Hence $R_m(z) = A_m(z)^{-1}$ is well defined on $F_{s+\epsilon,r-\delta}$.  Thus we have
\begin{equation}
\label{analyticres}
A(z)^{-1} = [A_m(z) + \rho_m(z)]^{-1} = [I + A_m(z)^{-1}\rho_m(z)]^{-1} A_m(z)^{-1}
\end{equation}
for all $z \in F_{s+\epsilon,r-\delta}$.

\section{Summary}
\label{sum}

We have established the key spectral properties of operator pencils using only elementary arguments based on the fundamental equations.  This insight may be important for applied mathematicians and engineers.   For matrix operators the fundamental equations can always be solved for $\{R_{-1}, R_0 \}$ near an isolated pole of the resolvent using Gaussian elimination \cite{alb1,how5}.  An interesting open question is whether they can always be solved for an infinite-dimensional operator near an isolated essential singularity of the resolvent.

\section*{Acknowledgements}
\label{ack}

This research was funded by ARC Discovery Grant \#DP1096551 held by Phil Howlett and Charles Pearce.  We thank an anonymous referee for several helpful suggestions.

\end{document}